\documentclass[12pt]{amsart}
 \usepackage[dvips]{graphics}
 \usepackage{amssymb, amsmath, mathrsfs}
 \input xy
 \xyoption{all}
 \DeclareMathAlphabet{\mathpzc}{OT1}{pzc}{m}{it}

 \newtheorem{theorem}{Theorem}[section]
 \newtheorem{proposition}[theorem]{Proposition}
 
 \newtheorem{lemma}[theorem]{Lemma}

 \newtheorem{conjecture}[theorem]{Conjecture}
 \theoremstyle{definition}
 \newtheorem{definition}[theorem]{Definition}

 \theoremstyle{remark}

 \newtheorem{remarks}[theorem]{Remarks}

\def\varle{\leqslant}

\newcommand{\CA}{{\mathcal A}}
\newcommand{\CB}{{\mathcal B}}

\newcommand{\CE}{{\mathcal E}}

\newcommand{\CG}{{\mathcal G}}
\newcommand{\CH}{{\mathcal H}}
\newcommand{\CI}{{\mathcal I}}

\newcommand{\CM}{{\mathcal M}}

\newcommand{\CO}{{\mathcal O}}

\newcommand{\CS}{{\mathcal S}}
\newcommand{\CT}{{\mathcal T}}

\newcommand{\CV}{{\mathcal V}}
\newcommand{\CW}{{\mathcal W}}

\newcommand{\CZ}{{\mathcal Z}}

\newcommand{\fh}{{{\mathfrak h}}}

\newcommand{\fg}{{{\mathfrak g}}}

\newcommand{\fn}{{{\mathfrak n}}}

\newcommand{\SA}{{\mathscr A}}
\newcommand{\SB}{{\mathscr B}}
\newcommand{\SC}{{\mathscr C}}

\newcommand{\SL}{{\mathscr L}}
\newcommand{\SM}{{\mathscr M}}
\newcommand{\SN}{{\mathscr N}}

\newcommand{\DR}{{\mathbb R}}
\newcommand{\DZ}{{\mathbb Z}}

\newcommand{\DN}{{\mathbb N}}

\newcommand{\Hom}{{\operatorname{Hom}}}
\newcommand{\id}{{\operatorname{id}}}

\newcommand{\udim}{{\operatorname{\underline{dim}}}}

\newcommand{\supp}{{\operatorname{supp}\,}}
\newcommand{\catmod}{{\operatorname{-mod}}}

\newcommand{\inj}{{\hookrightarrow}}
\newcommand{\Sh}{{\mathcal{SH}}}

\newcommand{\Loc}{{\SL}}

\newcommand{\linie}{{\,\text{---\!\!\!---}\,}}
\newcommand{\llinie}{{\text{---\!\!\!---\!\!\!---}}}

\begin{document}

\title{Kazhdan--Lusztig combinatorics via sheaves on Bruhat graphs}

\author{Peter Fiebig}
\address{Mathematisches Institut, Universit{\"a}t Freiburg, 79104 Freiburg, Germany}

\email{peter.fiebig@math.uni-freiburg.de}

%    General info
%\date{January 1, 1994 and, in revised form, June 22, 1994.}


\begin{abstract}
 To any Coxeter system we associate an exact category and  study the projective objects therein. We discuss an analog of the Kazhdan-Lusztig conjecture and show how it follows from a ``genericity'' conjecture and how the latter follows from a ``Hard Lefschetz'' conjecture.
\end{abstract}

\maketitle

\section{Introduction}
Recently essential progress in the theory of sheaves on polyhedral fans was made (i.e., the proof of the Hard Lefschetz Theorem by Karu and the combinatorial Koszul duality by Braden and Lunts). An analogous theory for Bruhat graphs of Coxeter systems still awaits a similar success. In this note we want to present some of the most intriguing open problems. 

 To a Coxeter system $(\CW,\CS)$ one  associates the {\em Hecke algebra} $\CH=\CH(\CW,\CS)$ over the ring $\DZ[v,v^{-1}]$ of Laurent polynomials in one variable. As a $\DZ[v,v^{-1}]$-module it is free with a basis $\{T_x\}$ indexed by elements in $\CW$. There is another basis, the  {\em Kazhdan-Lusztig self-dual basis} $\{C_x^\prime\}$. If $(\CW,\CS)$ is crystallographic, the  entries of the base change matrix encode graded multiplicities of representation theoretic and geometric objects. 

We present recent work on a combinatorial categorification of the above structures (\cite{Soe3, Fie2}). To $(\CW,\CS)$ we associate a {\em moment graph} $\CG=\CG(\CW)$ and give an algorithm that constructs a {\em sheaf} $\SB(x)$ on $\CG$ for any $x\in\CW$. We associate to $\SB(x)$ its {\em graded character} in the Hecke algebra and conjecture that it equals $C_x^\prime$.  This is equivalent to a conjecture that gives a bound on the degrees of the generators of the stalks of $\SB(x)$. It is well-known that, if $(\CW,\CS)$ is crystallographic,  the objects $\SB(x)$ correspond to intersection cohomology sheaves on the associated flag variety (\cite{Soe1, BMP}) and the conjecture then follows from work of Kazhdan and Lusztig (\cite{KL2}). In the non-crystallographic case the conjecture is very much open. The conjecture is proven in the dihedral cases in \cite{Soe3} and in the universal cases in \cite{Fie2}. 

Even in the crystallographic case there are intriguing problems left unanswered. Instead of a field of characteristic zero we might also use a field of positive characteristic as the ground field. (The characteristic should be such that the moment graph is a GKM-graph, cf.~ \cite{Fie1}.) In this case the conjecture corresponds to parts  of Lusztig's modular conjecture (\cite{Soe2}). 

For formal reasons one can deduce the conjecture for almost all characteristics from its characteristic zero analogue, though in general the exceptional primes are not known. One hopes that the conjecture is true for any characteristic that is at least the associated Coxeter number, but there is very little computational evidence yet.

In this note we will discuss the conjecture and show that it follows from a possibly stronger, Hard Lefschetz type conjecture. The latter is motivated by the Lefschetz condition in algebraic geometry. 
 
A proof of the conjecture would imply the positivity of the coefficients of the Kazhdan-Lusztig polynomials. Even in the crystallographic cases this could not yet be proven without geometric methods. Moreover, the conjecture implies that the Kazhdan-Lusztig polynomial $P_{y,x}$ itself only depends on the labelled Bruhat graph of the interval $[y,x]$ and thus supports a positive answer to a question of Dyer and Lusztig (cf.~ Section \ref{Dycon}).

\section{Sheaves on moment graphs}

\subsection{The moment graph associated to a Coxeter system}

Let $(\CW,\CS)$ be a Coxeter system, i.e.~ a group $\CW$ generated by a finite set $\CS\subset\CW$ with relations $s^2=1$ for each $s\in\CS$ and, possibly, additional relations of the form  $(st)^{m_{st}}=1$ with $m_{st}\geq 2$ and $s,t\in\CS$. Let $V$ be the {\em geometric representation}  of $(\CW,\CS)$ (\cite[5.3]{Hum}). It is a real, finite dimensional representation of $\CW$. Let $\CT\subset\CW$ be the set of {\em reflections} of $\CW$, i.e.~ the orbit of $\CS$ under the conjugation action. Then every $t\in\CT$ acts on $V$ as a reflection, i.e.~ it fixes a hyperplane pointwise and sends some non-zero vector $\alpha_t\in V$ to its negative. Note that $\alpha_t$ and $\alpha_s$ are linearly dependent only if $s=t$. We consider the contragredient representation on the dual vector space $V^\ast$. It is a faithful representation of $\CW$ such that every $t\in\CT$ acts as a reflection and, if we denote by $(V^\ast)^t$ the hyperplane stabilized pointwise by $t\in\CT$, then $(V^\ast)^s=(V^\ast)^t$ implies $s=t$. Note that $\alpha_t\in V$ is an equation of the reflection hyperplane $(V^\ast)^t$. It is well defined up to scalars.

To this datum we associate the labelled {\em Bruhat graph} $\CG=\CG(\CW,\CS)$ over $V$. It is one of the important examples of a {\em moment graph} (cf.~ \cite[1.2]{BMP} or \cite[1.2]{GKM} for an explanation of this terminology).  A moment graph is a directed graph without cycles whose edges are labelled by one-dimensional subspaces of $V$. In our situation it is constructed as follows. Let $\CW$ be its set of vertices  and connect $x$ and $y\in\CW$ by an edge $E$ if $x=ty$ for some $t\in\CT$. The edge is labelled by the line  $\DR\cdot \alpha_t\in V$ and is directed towards $y$ if $x<y$ in the Bruhat order. We denote it by $E\colon x\stackrel{\alpha_t}\to y$. 

\subsection{Sheaves on moment graphs}

Let $S$ be the algebra of real polynomial functions on $V^\ast$, i.e.~ the real symmetric algebra of the vector space $V$. There is an algebra $\DZ$-grading on $S$ that is determined by imposing that $V\subset S$ is the homogeneous part of degree $2$. We will only consider $\DZ$-graded $S$-modules $M$ in the following. Denote by $M_{\{k\}}$ the homogeneous part of degree $k$. Morphisms $f\colon M\to N$ of graded $S$-modules will always be of degree $0$, i.e.~ such that $f(M_{\{k\}})\subset N_{\{k\}}$ for all $k\in\DZ$.  

A {\em sheaf on $\CG$} is a combinatorial object $\SM=\{\SM^x,\SM^E,\rho_{x,E}\}$, where
\begin{enumerate}
\item $\SM^x$ is an $S$-module associated to any vertex $x$,
\item $\SM^E$ is an $S$-module associated to any edge $E$ with $\alpha\cdot \SM^E=0$ if $E$ is labelled by $\DR\cdot \alpha$,
\item for any vertex $x$ of $E$, $\rho_{x,E}\colon\SM^x\to\SM^E$ is a morphism of $S$-modules.
\end{enumerate}
In addition we assume that each $\SM^x$ is a torsion free and finitely generated $S$-module and that $\SM^x$ is zero for all but finitely many $x\in\CW$. A morphism $f\colon \SM\to\SN$ of sheaves is given by maps $f^x\colon \SM^x\to\SN^x$ and $f^E\colon\SM^E\to\SN^E$ of $S$-modules that are compatible with the maps $\rho$. Let $\Sh(\CG)$ be the resulting category of sheaves on $\CG$.

\section{IC-sheaves and indecomposable projective objects}
\subsection{The Braden-MacPherson construction}
Let $x\in\CW$. We recall the construction of the {\em intersection cohomology sheaf} $\SB(x)$. The terminology stems from the fact that, in the crystallographic case,  $\SB(x)$ encodes the structure of the equivariant intersection cohomology $IC_T(\overline{BxB/B})$ of the associated Schubert variety (\cite{BMP}). An analogous construction for rational fans in \cite{BBFK} gives the intersection cohomology of the associated toric variety. 

We set $\SB(x)^y=0$ if $y\not\leq x$ and $\SB(x)^x=S$. If $\SB(x)^y$ is already constructed for $y\in\CW$ and if $E\colon y^\prime\stackrel{\alpha}\to y$ is a directed  edge {\em ending} at $y$ (i.e.~ $y^\prime<y$), we set $\SB(x)^E=\SB(x)^y/ \alpha\cdot \SB(x)^y$ and we let $\rho_{y,E}\colon \SB(x)^y\to\SB(x)^E$ be the canonical quotient. We are left with constructing the $S$-modules $\SB(x)^y$ for $y<x$ and the maps $\rho_{y,E}$ for edges $E$ that originate at $y$. 

Let $y\in\CW$ and let $\CG_{>y}\subset\CG$ be the full subgraph with set of vertices $\{>y\}=\{w\in\CW\mid w>y\}$. Suppose that $\SB(x)$ is already  constructed on the subgraph $\CG_{>y}$. Define the {\em sections of $\SB(x)$} on $\CG_{>y}$ by
$$
\Gamma(\{>y\},\SB(x))=\left\{(m_w)\in\prod_{w>y}\SB(x)^w\left|\, \begin{matrix}\rho_{w,E}(m_w)=\rho_{w^\prime,E}(m_{w^\prime}) \\ \text{for all edges $E\colon w\linie w^\prime$ of $\CG_{>y}$}\end{matrix}\right.\right\}.
$$ 
Let $\CE^{\delta y}$ be the set of edges  originating at $y$ and let $\SB(x)^{\delta y}$ be the image of the map
$$
\Gamma(\{>y\},\SB(x))\stackrel{\bigoplus \rho_{w,E}}\longrightarrow \bigoplus_{E\in\CE^{\delta y}}\SB(x)^E.
$$
We define $\SB(x)^y$ as a projective cover of $\SB(x)^{\delta y}$ in the category of graded $S$-modules, and $\rho_{y,E}\colon \SB(x)^y\to\SB(x)^E$ as the components of the projective cover map $\SB(x)^y\to\SB(x)^{\delta y}$ for each edge $E\in\CE^{\delta y}$. This finishes the construction of $\SB(x)$.
\subsection{Sheaves with a Verma flag}
Next we want to define a category $\CV$ of sheaves on $\CG$ together with an exact structure. An exact structure allows us to talk about exact functors and {\em projective} objects in $\CV$. The IC-sheaves that we constructed in the last section will turn out to be the indecomposable projectives in $\CV$.

Let $\SM$ be any sheaf on $\CG$. The {\em space of global sections} of $\SM$ is
$$
\Gamma(\SM):=
\left\{ 
(m_x)\in\prod_{\text{$x$ a vertex of $\CG$}}  \SM^x\left| \,
\begin{matrix}
\rho_{x,E}(m_x)=\rho_{y,E}(m_{y}) \\
\text{for all edges $E\colon x\linie y$}
\end{matrix}
\right.
\right\}.
$$
The {\em structure algebra} of $\CG$ is defined by 
$$
\CZ=\CZ(\CG):=\left\{ 
(z_x)\in\prod_{\text{$x$ a vertex of $\CG$}}  S\left| \,
\begin{matrix}
z_x\equiv z_y\mod \alpha  \\
\text{for all edges $E\colon x\stackrel{\alpha}\llinie y$}
\end{matrix}
\right.
\right\}.
$$
The global sections of any sheaf $\SM$ on $\CG$ naturally form a $\CZ$-module by pointwise multiplication, since the label of each edge $E$ annihilates $\SM^E$.

There is a {\em localization functor} $\Loc$ that produces a sheaf on $\CG$ for any $\CZ$-module and that is left adjoint to $\Gamma$ (\cite{Fie1}). Let $\Loc\circ\Gamma\to \id$ be the adjunction morphism. We say that a sheaf $\SM$ {\em is generated by global sections} if $\Loc\circ\Gamma(\SM)\to\SM$ is an isomorphism.

Let $\CI\subset \CW$ be upwardly closed (i.e.~ each $x\in\CW$ which is bigger than some $y\in\CI$ in the Bruhat order also lies in $\CI$). Let $\SM$ be a sheaf on $\CG$ and let $\SM|_\CI$ be its restriction to the full subgraph $\CG_\CI$ with vertices in $\CI$. $\SM$ is said to be {\em flabby} if it is generated by global sections and if for any upwardly closed set $\CI\subset\CW$ the restriction $\Gamma(\SM)\to\Gamma(\SM|_\CI)$ is surjective. For a sheaf $\SM$ on $\CG$ and $w\in\CW$ let $\SM^{[w]}$ be the kernel of the restriction $\SM^w\stackrel{\bigoplus\rho_{w,E}}\longrightarrow\bigoplus_{E\in\CE^{\delta w}} \SM^E$.

\begin{definition} Let $\SM$ be a sheaf on $\CG$. We say that it {\em admits a Verma flag} if it is flabby and if, for any $w\in\CW$, the $S$-module $\SM^{[w]}$ is graded free. Let $\CV\subset\Sh(\CG)$ be the full subcategory of sheaves that admit a Verma flag. 
\end{definition}

\begin{theorem}[\cite{Fie2}] For any $x\in\CW$ the sheaf $\SB(x)$ admits a Verma flag.
\end{theorem}

The category $\CV$ is additive, but not abelian in all non-trivial cases. But it carries a canonical {\em exact structure} in the sense of Quillen. An exact structure is a class of sequences (called {\em short exact}) having similar properties as the class of short exact sequences in an abelian category (\cite{Qu}).

\begin{theorem}[\cite{Fie1}] There is an exact structure on $\CV$ consisting of those sequences $\SA\to\SB\to\SC$ in $\CV$ for which
$$
0\to\SA^{[w]}\to\SB^{[w]}\to\SC^{[w]}\to0
$$
is a short exact sequence of $S$-modules for any $w\in\CW$.
\end{theorem}

Let $\CA$ and $\CB$ be additive categories with exact structures and let $F\colon\CA\to\CB$ be a functor. It is called {\em exact} if it maps a short exact sequence to a short exact sequence. An object $P\in\CA$ is called {\em projective} if the functor $\Hom_{\CA}(P,\cdot)\colon\CA\to\DZ\catmod$ is an exact functor.

There are two characterizations of $\SB(x)$. The first one is local, while the second one is global. For $\SM\in\CV$ define $\supp\SM=\{x\in\CW\mid\SM^x\neq 0\}$. 
\begin{theorem}\begin{enumerate}
\item  $\SB(x)$ is the up to isomorphism unique indecomposable sheaf on $\CG$ with the following properties:
\begin{enumerate} 
\item For any $y\in\CW$, $\SB(x)^y$ is graded free over $S$,
\item $\supp\SB(x)\subset\{y\in\CW\mid y\leq x\}$ and $\SB(x)^x\cong S$, 
\item if $E\colon w\stackrel{\alpha}\to y$ is an edge, then $\SB(x)^E\cong\SB(x)^y/\alpha\cdot\SB(x)^y$ and $\SB(x)^w\to \SB(x)^{\delta w}$ is surjective.
\end{enumerate}
\item $\SB(x)\in\CV$ is the up to isomorphism unique object with the following properties:
\begin{enumerate}
\item $\SB(x)$ is indecomposable and projective.
\item $\supp\SB(x)\subset\{y\mid y\leq x\}$ and $\SB(x)^x\cong S$. 
\end{enumerate}
\end{enumerate}
\end{theorem}
A proof of the first claim can be found in  \cite{BMP} and a proof of the second in \cite{Fie1}.
\section{Connections to representation theory and geometry}
The following results are not needed in the sequel, but provide a motivation for the study of $\CV$. Let $\fg$ be a symmetrizable Kac--Moody algebra with associated Coxeter system $(\CW,\CS)$. The triangular decomposition $\fg=\fn^-\oplus \fh\oplus \fn$ of $\fg$ gives rise to the notion of highest weight module. The universal highest weight modules are called Verma modules. Let $\CM$ be the principal block of the category of $\fg$-modules that are extensions of finitely many Verma modules and that are semisimple for the action of $\fh$.

In \cite{Fie1} it is shown that $\CV$ is equivalent, as an exact category, to an equivariant version of $\CM$. Such an equivalence maps projective objects to projective objects, and hence the $\SB(x)$ provide a combinatorial model for the projective objects in $\CM$.

Suppose that  $\fg$ is finite dimensional and let $X$ be the flag variety associated to the Langlands dual Lie algebra. Braden and MacPherson showed that the sheaves $\SB(x)$ model the intersection cohomology sheaves of the Schubert varieties in $X$ (\cite{BMP}). 

Taken together, the two results above give the well-known correspondence between projective objects in $\CM$ and intersection cohomology sheaves on the associated Langlands dual flag variety.

\section{A degree conjecture and connections to Kazhdan-Lusztig combinatorics}

 For a graded $S$-module $M$ and $k\in\DZ$ let $M_{\{\varle k\}}$ be the submodule {\em generated in degrees $\leq k$}, i.e.~ $M_{\{\varle k\}}=\sum_{j\leq k} S\cdot M_{\{j\}}$. 

\begin{conjecture}[\cite{Fie2}]\label{dcon} $\SB(x)^y$ is generated in degrees $< l(x)-l(y)$, i.e.~
$$
\SB(x)^y=\SB(x)^{y}_{\{\varle l(x)-l(y)-1\}}.
$$
 \end{conjecture}

\begin{remarks}
\begin{enumerate}
\item For finite crystallographic Coxeter systems the conjecture follows from the Braden--MacPherson realization of intersection cohomology sheaves of Schubert varieties on the moment graph. In these cases the above inequality is translated to one of the defining axioms of intersection cohomology.
\item The constructions of $\CG$ and $\SB(x)$  make sense if we replace $V^\ast$ by any other reflection representation over an arbitrary field (cf.~ \cite{Fie2}).  If the field's characteristic is positive, the conjecture is related to Lusztig's modular conjecture, as explained in \cite{Soe2}.  In these cases the conjecture is wide open. However, for formal reasons the conjecture is true for almost all characteristics, provided it is true in characteristic zero.
\item The conjecture is true for dihedral Coxeter systems (\cite{Soe3}) and universal Coxeter systems (\cite{Fie2}), regardless of the ground field.
\end{enumerate}
\end{remarks}

Recently, a similar combinatorial theory associated to toric varieties instead of flag varieties was studied very successfully (cf.~ \cite{B} for a review). In particular, the analog of Conjecture \ref{dcon} was proven by Karu in \cite{Ka} without the use of geometry.

The next section shows how one can relate the conjecture to combinatorial problems in the Hecke algebra associated to $(\CW,\CS)$. 

\subsection{Kazhdan-Lusztig polynomials}

Let $\CH$ be the Hecke algebra of $(\CW,\CS)$, i.e.~ the free $\DZ[v,v^{-1}]$-module with basis elements $T_x$ for $x\in\CW$, equipped with a multiplication such that
\begin{eqnarray*}
{T}_x\cdot {T}_{y} & = & {T}_{xy}\quad\text{if $l(xy)=l(x)+l(y)$}, \\
{T}_s^2 & = & v^{-2}T_e+(v^{-2}-1)T_s \quad\text{for $s\in\CS$}.
\end{eqnarray*}
Note that $T_e$ is a unit in $\CH$. One can show that there exists an inverse $T_w^{-1}$ for any $w\in\CW$. There is an involution $d\colon \CH\to\CH$ given by $d(v)=v^{-1}$ and $d(T_w)=T_{w^{-1}}^{-1}$ for $w\in\CW$. Set $\tilde T_x=v^{l(x)}T_x$.

\begin{theorem}[\cite{KL1}]\label{KLelts} For any $x\in\CW$ there exists a unique element $C_x^\prime=\sum_{y\in\CW} h_{y,x}(v)\cdot \tilde T_y\in\CH$ with the following properties:
\begin{enumerate}
\item \label{KLeltsd} $C_x^\prime$ is self-dual, i.e.~ $d(C^\prime_x)=C^\prime_x$. 
\item \label{KLeltss} $h_{y,x}(v)=0$ if $y\not\leq x$, and $h_{x,x}(v)=1$,
\item \label{KLeltsn} $h_{y,x}(v)\in v\DZ[v]$ for $y<x$. 
\end{enumerate}
\end{theorem} 

There are polynomials $P_{y,x}\in\DZ[v,v^{-1}]$ such that $P_{y,x}(v^{-2})=v^{l(y)-l(x)}h_{y,x}(v)$. They are called the Kazhdan--Lusztig polynomials.

\subsection{Characters in $\CH$}
 
For $l\in\DZ$ the {\em shift functor} $M\mapsto M\{l\}$ on graded $S$-modules shifts the grading down by $l$, i.e.~ such that $M\{l\}_{\{k\}}=M_{\{l+k\}}$ for $k\in\DZ$. Let $x,y\in\CW$ and suppose that $\SB(x)^{[y]}\cong\bigoplus_i S\{k_i\}$ for $k_1,\dots,k_n\in\DZ$. The multiset of elements $\{k_1,\dots,k_n\}$ is well-defined, so we can define the {\em graded rank} of $\SB(x)^{[y]}$  by $\udim(\SB(x)^{[y]})=\sum_{i=1}^n v^{-k_i}\in\DZ[v,v^{-1}]$. We define the {\em graded character} of $\SB(x)$ by
$$
h(\SB(x))=\sum_{y\in\CW}\udim(\SB(x)^{[y]})\cdot v^{l(y)-l(x)}\tilde T_y\in\CH.
$$

The following conjecture is a variant of Kazhdan and Lusztig's conjecture in \cite{KL1}. It also appears in \cite{Soe3} in the context of ``$S$-bimodules on twisted diagonals''.
 
\begin{conjecture}\label{KLcon} We have $h(\SB(x))=C_x^\prime$ for any $x\in\CW$. 
\end{conjecture}

There is a covariant equivalence
$D\colon \CV^\circ\cong\CV^{opp}$,
where $\CV^{\circ}\subset\CZ\catmod^f$ is the full subcategory of modules that admit a Verma flag with respect to the {\em dual} partial ordering on $\CW$, and where $\CV^{opp}$ denotes the opposite category of $\CV$. On the level of global sections it is given by the duality $(\cdot)^\ast\colon\CZ\catmod^f\to\CZ\catmod^f$, $M^\ast=\bigoplus_i \Hom_S(M,S\{i\})$, with the $\CZ$-action defined by $(z.f)(m)=f(z.m)$ for $z\in\CZ$, $f\in M^\ast$ and $m\in M$. 

\begin{theorem}\label{duality} Choose $x\in\CW$.
\begin{enumerate}
\item There is an isomorphism $D(\SB(x))\cong \SB(x)\{2l(x)\}$, i.e.~ $\SB(x)$ is self--dual of degree $2l(x)$.
\item The isomorphism above induces an isomorphism 
$$
\SB(x)^{[y]}\cong D(\SB(x)^y)\{-2l(x)+2l(y)\}.
$$
\end{enumerate}
\end{theorem}
In \cite{Fie2} this is proven using translation functors. Hence conjecture \ref{dcon} is true if and only if  $\SB(x)^{[y]}$ lives in degrees $>l(x)-l(y)$.

\begin{proposition}[\cite{Fie2}] Conjectures \ref{dcon} and \ref{KLcon} are equivalent.
\end{proposition}
\begin{proof}[Sketch of a proof] Proving Conjecture \ref{KLcon} amounts to showing that the character $h(\SB(x))$ has the defining properties of $C^\prime_x$ that are listed in Theorem \ref{KLelts}. 
The self-duality of Theorem \ref{duality} implies that $h(\SB(x))$ is a self-dual element in $\CH$. It is clear from the construction of $\SB(x)$ that $h(\SB(x))$ is supported on $\{\leq x\}$. And the normalization property (\ref{KLeltsn}) is just a reformulation of the statement of Conjecture \ref{dcon}.
\end{proof}

\subsection{Does $P_{y,x}$ depend only on the isomorphism type of the poset $[y,x]$?}\label{Dycon}

The combinatorial invariance conjecture, proposed independently by Dyer (\cite{Dye}) and Lusztig, reads:  
$P_{x,y}$ depends only on the isomorphism type of the Bruhat interval $[y,x]$. 

In \cite{BCM} Brenti, Caselli and Marietti gave strong support for the conjecture, generalizing  earlier results of du Cloux (\cite{dC}) and Brenti (\cite{Br2}), by showing that for $x,y\in\CW$, $y\leq x$, the polynomial $P_{y,x}$ depends only on the poset $[e,x]$.

In the crystallographic cases the conjecture \ref{KLcon} is proven and we deduce the following:

\begin{theorem}[{\cite[Theorem 7.3]{Br1}}]\label{pol} $P_{y,x}$ depends only on the Bruhat graph $\CG_{[y,x]}$.
\end{theorem}

In a sense this complements the combinatorial results in \cite{BCM} and \cite{dC}. While in both papers  the authors use only the poset structure of the interval $[e,x]$, Theorem \ref{pol} uses additional information on the moment graph, namely the labelling of the edges, but only on the interval $[y,x]$.

\section{A  genericity conjecture and a Hard Lefschetz conjecture}

Consider for any $x,y\in\CW$ with $y<x$ the inclusion $\SB(x)^{[y]}\inj\SB(x)^y$ of graded free $S$-modules.  

\subsection{Genericity}
Let $l\in\DZ$ and 
suppose that $M$ is a graded free $S$-module of multirank $\{-l-k_1,\dots,-l-k_n\}$ for some $k_1,\dots,k_n\geq 0$ (i.e.~ $M\cong \bigoplus S\{-l-k_i\}$) and that $N$ is a graded free $S$-module of multirank $\{-l+k_1,\dots, -l+k_n\}$. We want to consider all maps $f\colon M\to N$ of graded $S$-modules. Note that the space of such maps is a finite dimensional vector space. 

\begin{lemma}\label{genmaps} Generically, i.e.~for $f$ in a Zariski open non-empty set, $f$ is injective and we have
$$
f(M_{\{\varle  l+ m-1\}})\cap N_{\{\varle  l-m\}}=0.
$$
for all $m\geq 1$.
\end{lemma}

\begin{proof} Fix $m\geq 1$ and consider the canonical map 
$$
\rho_m\colon \Hom(M,N)\to \Hom(M_{\{\varle  l+m-1\}}, N/N_{\{\varle  l-m\}}).
$$ 
We have to show that $\rho_m(f)$ is injective for almost all $f$ and all $m\in\DN$.  For fixed $m$, $\rho_m$ is surjective and almost all maps in $\Hom(M_{\{\varle  l+m-1\}}, N/N_{\{\varle  l-m\}})$ are injective. For large $m$, $M_{\{\varle  l+m-1\}}=M$ and $N/N_{\{\varle  l-m\}}=N$, hence we only have to consider finitely many $m$. The claim follows.
\end{proof}

\subsection{An inductive assumption}
Let $x\in\CW$.
In \cite{Fie2} it is shown how to approach Conjecture \ref{dcon} inductively on the Bruhat order. In particular it follows from the validity of the conjecture for  some $x^\prime$ with $x^\prime<x$ and $x^\prime=xs$ for $s\in\CS$, that $\SB(x)^y$ is generated in degrees $\leq l(x)-l(y)$. From now on we will take this as an assumption.

Choose $y\in\CW$ with $y<x$. Then $\SB(x)^y\cong\bigoplus S\{-(l(x)-l(y))+k_i\}$ for some $k_i\geq 0$. From the duality in \ref{duality} it follows that $\SB(x)^{[y]}=\bigoplus S\{-(l(x)-l(y))-k_i\}$, so we can apply the preceding discussion to $M=\SB(x)^{[y]}$ and $N=\SB(x)^y$ with $l=l(x)-l(y)$. We conjecture

\begin{conjecture}\label{pcon} For any $x,y\in\CW$ with $y<x$ the inclusion $\SB(x)^{[y]}\to \SB(x)^y$ is generic in the sense of Lemma \ref{genmaps}, i.e.~ for any $m\geq 1$ we have
$$
\SB(x)^{[y]}_{\{\varle l(x)-l(y)+m-1\}}\cap \SB(x)^y_{\{\varle l(x)-l(y)-m\}}=0.
$$
\end{conjecture}

\begin{proposition} Conjecture \ref{pcon} implies Conjecture \ref{dcon}.
\end{proposition}
\begin{proof} We have to show that $\SB(x)^y=\SB(x)^y_{\{\varle l(x)-l(y)-1\}}$. By duality this is equivalent to $\SB(x)^{[y]}_{\{\varle l(x)-l(y)\}}=0$. Assume that Conjecture \ref{pcon} is true. Then, in particular, the map
$$
\SB(x)^{[y]}_{\{\varle l(x)-l(y)\}}\to\SB(x)^y/\SB(x)^y_{\{\varle l(x)-l(y)-1\}}
$$
is injective. Suppose $v\in \SB(x)^{[y]}_{\{l(x)-l(y)\}}$ is non-zero. Then the image of $v$ in $\SB(x)^y$ generates a direct summand isomorphic to $S\{-(l(x)-l(y))\}$. This is impossible, since $v$ maps to zero under the projective cover map $\SB(x)^y\to \SB(x)^{\delta y}$.
\end{proof}

\subsection{A Hard Lefschetz conjecture}

Choose a line $\DR\cdot v\subset V^\ast$ that is not contained in any reflection hyperplane. It correponds to a surjective map $\rho\colon S\to \DR[T]$ such that $\rho(\alpha_t)\neq 0$ for all $t\in\CT$. This is a map of graded algebras if we give $\DR[T]$ the grading such that $\DR\cdot T$ is the homogeneous part of degree $2$.  For any $S$-module $M$ let $M_T:=M\otimes_S \DR[T]$ be its restriction to the chosen line.  For $v\in M$  denote by $\bar v=v\otimes 1\in M_T$ its image. 

Let $\alpha_1,\dots,\alpha_k\in V$ be the labels of all edges that originate at $y$ in the graph $\CG_{\leq x}$. Then multiplication with $\alpha_1\cdots\alpha_k\in S$ maps $\SB(x)^y$ into $\SB(x)^{[y]}$. Hence multiplication with $T^k$ maps $\SB(x)^y_T$ into the image of the map $\SB(x)^{[y]}_T\to\SB(x)^y_T$. So $\SB(x)^{\delta y}_T$ is a torsion $\DR[T]$-module. Since $\SB(x)^{y}$ and $\SB(x)^{[y]}$ have the same ungraded rank, we deduce that $\SB(x)^{[y]}_T\to\SB(x)^y_T$ is injective. In particular, if $v\in\SB(x)^{[y]}$ is part of a graded $S$-basis, then the image of $v$ in $\SB(x)^y_T$ is non-zero. We conjecture

\begin{conjecture}[cf.~ {\cite[Bemerkung 7.2]{Soe3}}]\label{p-dimone} For any $x,y\in\CW$ with $y<x$ the $\DR[T]$-module $\SB(x)^{\delta y}_T$ is isomorphic  to a direct sum of modules of the form $(\DR[T]/T^{n+1})\{-(l(x)-l(y))+n\}$ for $n\geq 0$. 
\end{conjecture}

\begin{proposition} Conjecture \ref{p-dimone} implies Conjecture \ref{pcon}.
\end{proposition}
\begin{proof} Fix  $m\geq 1$ and suppose  
$$
v\in \SB(x)^{[y]}_{\{\varle  l(x)-l(y)+m-1\}}\cap \SB(x)^y_{\{\varle  l(x)-l(y)-m\}}
$$ 
is non-zero. We can assume that  $v$ is homogeneous of degree $k\leq l(x)-l(y)+m-1$ and generates a direct summand in $\SB(x)^{[y]}$ that is isomorphic to $S\{-k\}$. Then there are $v_1,\dots, v_l\in\SB(x)^y_{\{l(x)-l(y)-m\}}$ and $f_1,\dots,f_l\in S$, homogeneous of degree $k-(l(x)-l(y)-m)\leq 2m-1$, such that $x=f_1\cdot v_1+\dots+f_l\cdot v_l$. Specializing gives an equation 
$$
\bar x=\bar f_1\cdot \bar v_1+\dots+\bar f_l\cdot \bar v_l
$$
in $\SB(x)^{y}_T$. Now $\bar x\neq 0$ in $\SB(x)^y_T$, but its image in $\SB(x)^{\delta y}_T$ is zero. If we assume Conjecture \ref{p-dimone}, this contradicts the fact that $T^m$  acts injectively when restricted to $\SB(x)^{\delta y}_{T,\{\varle  l(x)-l(y)-m\}}$. 
\end{proof}

Multiplication by $T$ is a nilpotent operation on $\SB(x)_T^{\delta y}$ of degree $2$. We conjecture

\begin{conjecture}[Hard Lefschetz property]\label{HLcon} Let $x,y\in\CW$ with $y<x$. Multiplication by $T$ induces an isomorphism
$$
T^{m}\colon \SB(x)^{\delta y}_{T,\{l(x)-l(y)-m\}}\stackrel{\sim}\to\SB(x)^{\delta y}_{T,\{l(x)-l(y)+m\}}.
$$
for all $m\geq 1$.
\end{conjecture}

\begin{proposition} Conjecture \ref{HLcon} is equivalent to Conjecture \ref{p-dimone}.
\end{proposition}

\begin{proof} This is a  direct consequence of the following lemma.
\end{proof}

\begin{lemma} Let $M=\bigoplus_{i\in I} (\DR[T]/T^{n_i+1})\{-l_i\}$ be a graded $\DR[T]$-module and suppose that there is $l\in\DZ$ such that, for each $m\geq 1$,
$$
T^m\colon M_{\{l-m\}}\stackrel{\sim}\to M_{\{l+m\}}\eqno{(\ast)}
$$ 
is an isomorphism. Then $l_i=l-n_i$ for each $i\in I$.
\end{lemma}
\begin{proof} It is easily seen that $(\ast)$ holds if and only if it holds on every direct summand $(\DR[T]/T^{n_i+1})\{-l_i\}$. On such a direct summand $(\ast)$ holds if and only if $l_i=l-n_i$.
\end{proof}

\end{document}